\documentclass[a4paper,fleqn]{cas-dc}
\usepackage{amsmath,amsfonts}
\usepackage[ruled,linesnumbered]{algorithm2e}
\usepackage{multirow}
\usepackage[authoryear]{natbib}
\def\tsc#1{\csdef{#1}{\textsc{\lowercase{#1}}\xspace}}
\tsc{WGM}
\tsc{QE}

\begin{document}
\let\WriteBookmarks\relax
\def\floatpagepagefraction{1}
\def\textpagefraction{.001}

\shorttitle{}

\shortauthors{Yajie Wen}

\title [mode = title]{A Beam Search Based Parallel Algorithm for the Two-Dimensional Strip Packing Problem}

\author[1]{Yajie Wen}[
	orcid = 0009-0004-2053-0583,
]
\fnmark[1]
\ead{23020201153806@stu.xmu.edu.cn}

\author[2]{Defu Zhang}[
]
\cormark[1]
\fnmark[2]
\cortext[2]{Corresponding author}
\ead{dfzhang@xmu.edu.cn}

\affiliation{organization={Xiamen University},
            addressline={Siming south street 443}, 
            city={Xiamen}, 
            postcode={361000}, 
            state={Fujian},
            country={China}}

\begin{abstract}
This paper introduces BSPA, a parallel algorithm that leverages beam search to address the two-dimensional strip packing problem. The study begins with a comprehensive review of existing approaches and methodologies, followed by a detailed presentation of the BSPA algorithm. Experimental results demonstrate the effectiveness of the proposed method. To facilitate further research, both the code and datasets are publicly available.
\end{abstract}



\maketitle

\section{Introduction}
Given a set of \( n \) rectangular boxes and a container with a fixed width \( W \) and infinite length, the two-dimensional strip packing (2DSP) problem involves arranging all the boxes within the container such that the total length used by the boxes is minimized.

The two-dimensional strip packing problem is widely applied in industries such as metal processing, wood processing, glass manufacturing, and furniture production, where it is necessary to cut various shapes and sizes of small parts from large sheets of raw material. The 2DSP problem helps optimize cutting plans, enhancing raw material utilization and minimizing waste. For example, different components such as table tops and chair legs must be cut from standard-sized wood boards in furniture manufacturing. By applying an effective two-dimensional strip packing algorithm, the optimal cutting sequence and layout can be determined.

The 2DSP problem can be classified into four categories based on whether box rotation is allowed and whether guillotine cuts are required. Moreover, guillotine refers to a type of cut where the cutting action spans from one side to the other:
\begin{itemize}
	\item {\texttt{RF}}: Rotation is allowed, and no guillotine cut is required.
	\item {\texttt{RG}}: Rotation is allowed, and a guillotine cut is required.
	\item {\texttt{OF}}: Rotation is forbidden, and no guillotine cut is required.
	\item {\texttt{OG}}: Rotation is forbidden, and a guillotine cut is required.
\end{itemize}

This paper explicitly addresses the RF and OF cases, which require a guillotine cut, and considers the constraints on box rotation.

\subsection{Related work}
The two-dimensional strip packing problem (2DSP) has garnered significant attention over several decades due to its theoretical complexity (NP-hard) and practical importance in logistics, manufacturing, and resource optimization. Existing approaches to solving the 2DSP can be categorized into three primary paradigms: exact algorithms, heuristic algorithms, and emerging machine learning-based methods. The subsequent section offers a comprehensive overview of key advancements within each category.
\subsection{Exact algorithms}
Early research on the two-dimensional strip packing problem (2DSP) primarily focused on exact algorithms to guarantee optimal solutions, albeit at significant computational cost. Over time, these exact approaches were refined. \cite{Hifi1998} combined branch-and-bound with dynamic programming, while \cite{Martello2003} and \cite{Kenmochi2009} developed specialized branch-and-bound algorithms tailored specifically for the 2DSP. \cite{Cote2014} introduced column generation to obtain exact solutions. To address guillotine constraints, \cite{Messaoud2008} proposed polynomial-time feasibility checks, and \cite{Fleszar2016} further advanced the concept of stage-unrestricted guillotine cutting.

Efforts to reduce computational costs through parallelization were also explored. \cite{Bak2011} proposed a parallel branch-and-bound algorithm, and \cite{Boschetti2010} developed reduction procedures aimed at enhancing efficiency. Despite these advancements, exact methods remain impractical for large-scale instances due to their exponential time complexity.

While exact algorithms provide theoretical rigour, they primarily apply to small-scale or highly constrained instances, as their computational demands become prohibitively high for larger problem sizes.

\subsection{Heuristic algorithms}
Heuristic methods have become the predominant approach for solving the two-dimensional strip packing problem, offering a practical balance between computational efficiency and solution quality. Early contributions by \cite{Baker1980} introduced the Bottom-Up Left-Justified (BL) algorithm, which laid the foundation for approximation bounds in packing problems. This seminal work was further refined by \cite{Chazelle1983}, who optimized the bottom-left heuristic to enhance the efficiency of placement reporting. Building upon the BL algorithm, subsequent studies introduced several noteworthy variants. \cite{Liu1999} improved packing density by leveraging permutation-based patterns, while \cite{Wei2011} incorporated greedy selection techniques combined with tabu search, demonstrating the versatility of heuristics in this domain.

A significant advancement in heuristic methods was made by \cite{Burke2004a}, who proposed a best-fit heuristic for strip packing. This approach was subsequently extended by \cite{Asik2009} and \cite{Oezcan2013}, who introduced bidirectional niche placement and compound polygon evaluation, further improving the quality of the solutions. The application of genetic algorithms (GAs) as a metaheuristic gained prominence, with \cite{Kroeger1995} pioneering their use for rectangle packing problems. \cite{Yeung2004} later hybridized GAs with a left-fit-bottom heuristic, and \cite{Bortfeldt2006} refined this hybrid approach specifically for strip packing, showcasing the potential of evolutionary algorithms in solving complex packing problems. Iterative heuristics, such as those proposed by \cite{Belov2008} with the SVC and BLR methods and \cite{He2013} with a deterministic heuristic (DHA), have also contributed to the advancement of strip packing techniques.

In parallel, recursive and block-based strategies emerged as effective alternatives. \cite{Cui2008} developed the HRBB algorithm, applying recursive branch-and-bound techniques to the strip packing problem. Building on this, \cite{Cui2013} introduced the Sequential Grouping and Value Correction Procedure (SGVCP), which iteratively refines packing layouts to optimize overall efficiency. \cite{Wei2016} proposed a three-stage intelligent search algorithm that combines greedy, local, and randomized improvement methods, achieving superior performance compared to existing metaheuristic approaches.

The exploration of greedy randomized adaptive search (GRASP) methods also gained traction. \cite{AlvarezValdes2008} introduced a reactive GRASP algorithm, which was later refined by \cite{OviedoSalas2022}, further enhancing the performance of these techniques. Hybrid approaches, such as the Hybrid Heuristic Algorithm (HHA) proposed by \cite{Chen2019}, which integrates greedy and local search methods, and adaptive learning enhancements by \cite{Rakotonirainy2020}, have further advanced the effectiveness of metaheuristics in solving strip packing problems.

Recent innovations in the field include \cite{Zhang2024}, who introduced the Block-Based Heuristic Search Algorithm (BBHSA), optimizing block-level layouts and reducing fragmentation. \cite{Grandcolas2022} proposed the PVS strategy, merging local search with exact verification techniques. Additionally, \cite{Wei2017, Wei2019} developed the BestFitPack and FirstFitPack algorithms, focusing on bottom-left placement strategies and fitness evaluation to improve packing efficiency.

Together, these studies underscore the continuous development of heuristic and metaheuristic approaches for the two-dimensional strip packing problem. By leveraging problem-specific insights and algorithmic adaptations, these methods have achieved near-optimal solutions in a computationally efficient manner, highlighting the flexibility and effectiveness of heuristics in addressing this challenging optimization problem.

\subsection{Machine learning and hybrid approaches}
Recent advancements have increasingly focused on leveraging learning-based techniques to address the combinatorial complexity of the two-dimensional strip packing problem. \cite{Zhao} applied Q-learning to optimize packing sequences, improving material utilization. \cite{Fang2023} integrated deep reinforcement learning (DRL) with pointer networks, combining Maxrects-BL positioning with sequence optimization to enhance packing efficiency. \cite{NeuenfeldtJunior2021} introduced a meta-learning approach to select optimal multi-label transformations, further advancing the application of machine learning in solving 2DSP. These studies demonstrate the potential of learning-based methods to tackle the inherent challenges of the problem, offering new avenues for improving solution quality and computational efficiency.

The remainder of the paper is organized as follows: Section 2 presents the proposed beam-search-based approach. Section 3 details the computational experiments conducted to evaluate the algorithm's performance. Finally, Section 4 provides conclusions.

\section{Algorithm description}
Our algorithm is inspired by the work of \cite{Zhang2024}, which reformulates the two-dimensional strip packing problem as a two-dimensional rectangular packing problem (2DRP) with a fixed container size. Our approach decomposes 2DSP into a 2DRP problem and a smaller 2DSP subproblem, employing three beam searches to obtain the final solution. The first beam search packs boxes into a container with a length determined by the total area of all boxes divided by the container width, aiming to obtain an optimal packing configuration directly. The second beam search addresses the remaining unpacked boxes from the first stage, effectively solving a reduced 2DSP problem. The third beam search refines the solution by exploring container lengths within the range from the first container’s length to the sum of the lengths of the first and second containers.
Let \( W \) denote the container's width and \( S \) represent the total area of all boxes. The method proposed in this study consists of three main steps: block generation, minimum container filling, and multiple container filling.

\emph{Block generation:} In this step, rectangular boxes are systematically paired and combined to form larger rectangular blocks.

\emph{Minimum container filling:} A container \( C_a \) is generated with dimensions length \( S/W \) and width \( W \). The blocks formed in the previous stage are then used to fill container \( C_a \) using the beam search algorithm. The optimal filling configuration, denoted as \( bestState_a \), is recorded during this process.

\emph{Multiple container filling:} If all boxes are successfully packed into \( C_a \) according to \( bestState_a \), then \( bestState_a \) is considered the optimal solution. Otherwise, a new container \( C_b \) is generated with a width of \( W \) and infinite length. The remaining boxes from the second step are then used to fill container \( C_b \) via the beam search algorithm, and the best packing configuration is recorded as \( bestState_b \). The filling lengths of the containers in \( bestState_a \) and \( bestState_b \) are summed to obtain the maximum length, denoted as \( L_{max} \). Subsequently, a series of container lengths \( L_i \) are generated within the range \( \left( S/W, L_{max} \right] \), with a step size of 1. For each \( L_i \), a container \( C_i \) with dimensions \( L_i \times W \) is created, where \( i \) ranges from 1 to \( \lceil L_{max} - S/W \rceil \). Finally, the blocks are used to fill each container \( C_i \) using the beam search algorithm, and the best packing results are returned after the search is completed. Parallel computing can be used to search for \( C_i \).
\subsection{Spatial representation}\label{AA}
\begin{figure}[h]
	\centering
	\includegraphics[width=0.7\linewidth]{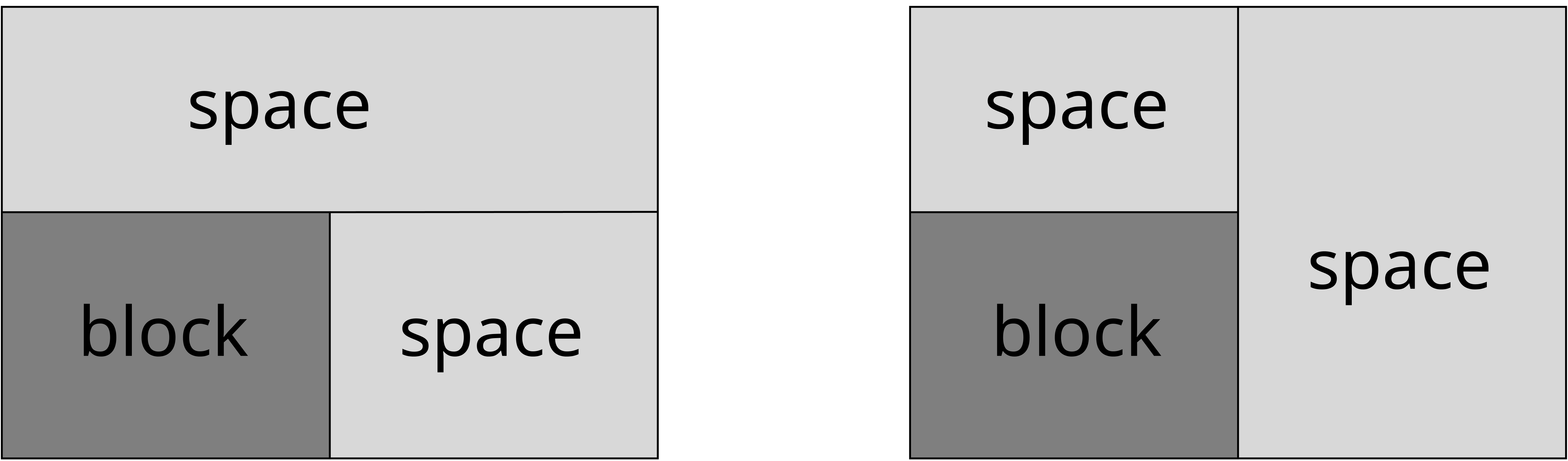}
	\caption{Partial Representation.}
	\label{pspace}
\end{figure}

\begin{figure}[h]
	\centering
	\includegraphics[width=0.7\linewidth]{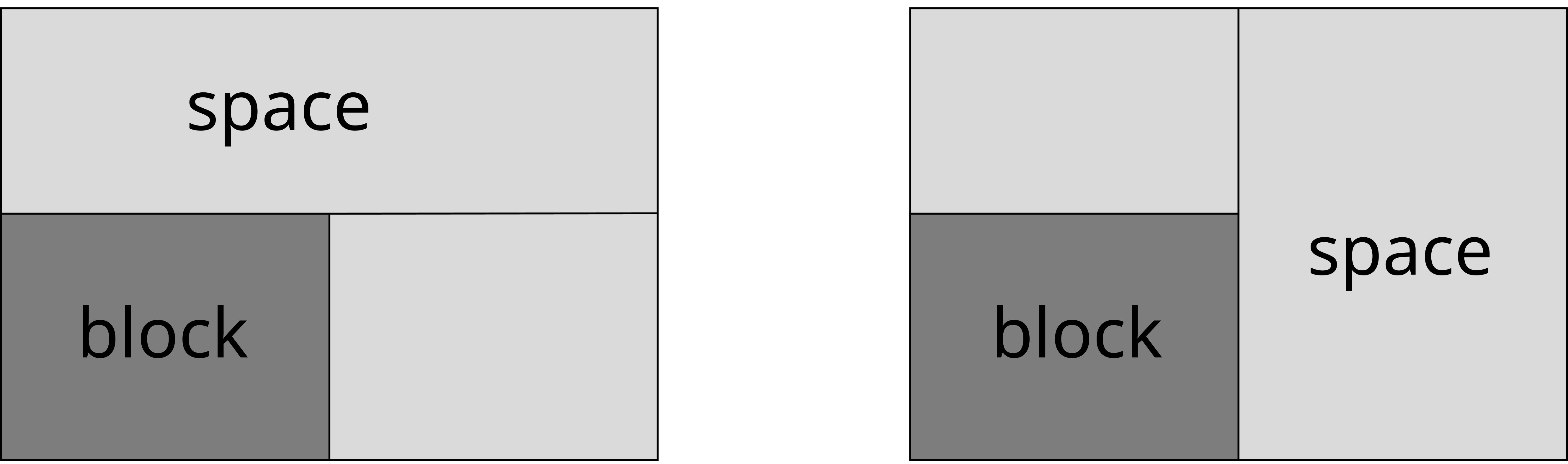}
	\caption{Overlapping Representation.}
	\label{cspace}
\end{figure}

Before presenting the algorithm, it is necessary to introduce two standard space representation methods: partial and overlapping representations. In partial representation, the spaces are non-overlapping, whereas in overlapping representation, the spaces can intersect. Figure \ref{pspace} illustrates the partial representation, while Figure  \ref{cspace} shows the overlapping representation. This study adopts the overlapping representation method, with further details on the space update process available in the paper \cite{Neveu2007}.

\subsection{Block generation}\label{AA}
The block generation approach employed in this study consists of two steps: the simple block generation method, as described in Algorithm \ref{sblock}, and the complex block generation method, outlined in Algorithm \ref{cblock}. A simple block is defined as a block composed of boxes of identical size, while a complex block is formed by combining multiple simple blocks. The process of creating complex blocks involves joining two simple blocks both lengthwise and widthwise, subject to the following conditions:
\begin{itemize}
	\item The length and width of the generated block must be smaller than the container's.
	\item The total number of boxes in the block must not exceed the available boxes.
	\item The ratio of the area occupied by the boxes in the block to the total area of the block must be at least \(minFillRate\).
	\item The total number of blocks must not exceed \(maxNum\).
\end{itemize}

For the block generation process in the study, the parameters are set as follows: \(maxNum = 10000\) and \(minFillRate = 1\).

\begin{algorithm}
	\caption{simpleBlocks}
	\label{sblock}
	\KwIn{$boxList$}
	\KwOut{$blocks$}
	
	$blocks \gets \emptyset$
	
	\For{each type of \(box\) in $boxList$}{
		$boxNum \gets$ The number of boxes of type $box$\;
		\For{$x \gets 1$ to $boxNum$}{
			\For{$y \gets 1$ to $\lfloor boxNum / x \rfloor$}{
				$block \gets$ Generates a block of width $x*boxWidth$ and length $y*boxLength$\;
				Add $block$ to $blocks$\;
			}
		}
	}
	Rotate the block in $blockList$ under the RF case \;
	\Return{$blocks$} \;
\end{algorithm}

\begin{algorithm}
	\caption{complexBlocks}
	\label{cblock}
	\KwIn{$boxList$, $maxNum$, $minFillRate$}
	\KwOut{$blockList$}
	
	$blockList \gets$ simpleBlocks($boxList$) \;
	$pList \gets blockList$ \;
	$newBlockList \gets \emptyset$ \;
	
	\While{$|blockList| < maxNum$}{
		$newBlockList \gets \emptyset$ \;
		\ForEach{$block_1 \in pList$}{
			\ForEach{$block_2 \in blockList$}{
				$combinedBlocks \gets$ Concatenate $block_1$ and $block_2$ in length and width directions according to the conditions\;
				\If{$combinedBlocks$ is not empty}{
					Add all blocks from $combinedBlocks$ to $newBlockList$\;
				}
			}
		}
		
		\If{$newBlockList$ is empty}{
			\textbf{break}
		}
		
		\ForEach{$block \in newBlockList$}{
			\If{$|blockList| < maxNum$}{
				Add $block$ to $blockList$\;
			}
		}
		
		$pList \gets newBlockList$ \;
	}
	\Return{$blockList$} \;
\end{algorithm}

\subsection{Minimum container filling}
Before initiating the minimum container filling process, a runtime for beam search is set, and an initial state, \( stateInit \), is generated. This initial state includes a list of remaining spaces, \( spaceList \), a list of blocks, \( blockList \), the number of available boxes, \( avaibox \), and a state score. Initially, the spaces in \( spaceList \) correspond to the entire container \( C_a \), \( avaibox \) represents the number of boxes of each size, and \( blockList \) contains the generated complex blocks.

In the beam search process of this step, starting from \( stateInit \), a space is selected from \( spaceList \) using the space selection method. Then, \( w \) blocks are chosen from \( blockList \) using the block selection method. These \( w \) blocks are placed into the lower-left corner of the selected space, resulting in \( w \) new states. The scores of these new states are computed using a greedy approach, and the best solution, denoted as \( bestState_a \), is recorded. Finally, the top \( w \) states with the highest scores are selected as the nodes for further expansion. The pseudocode for this process is outlined in Algorithm \ref{search1}.

\begin{algorithm}[ht]
	\caption{minSearch}
	\label{search1}
	\KwIn{$stateInit$, $bestState_a$}
	\KwOut{$bestState_a$}
	
	$w \gets 1$ \;
	
	\While{Time Limit Is Reached}{
		$stateList \gets \emptyset$ \;
		Add $stateInit$ to $stateList$ \;
		
		\While{$stateList$ is not empty}{
			$successors \gets \emptyset$ \;
			\ForEach{$state \in stateList$}{
				\If{$state$ equals $stateInit$}{
					$succ \gets$ expand($state$, $w \times w$) \;
				}
				\Else{
					$succ \gets$ expand($state$, $w$) \;
				}
				Add all states in $succ$ to $successors$ \;
			}
			
			$stateList \gets \emptyset$ \;
			
			\If{$successors$ is not empty}{
				\ForEach{$state \in successors$}{
					greedy($state$, $bestState_a$) \;
				}
			}
			Add the top $w$ states in $successors$ ranked by score to $stateList$ \;
		}
		
		$w \gets \sqrt{2} \times w$ \;
	}
	
	\Return{$bestState_a$} \;
\end{algorithm}
Next, the detailed specifics of the algorithm are presented.

\emph{Space selection method:} Let the coordinate origin be \( (0,0) \), and the lower-left corner of the space be \( (x, y) \). The space to be selected is the one that minimizes the sum \( x + y \), i.e., the space closest to the origin.

\emph{Block selection method:} First, identify all blocks in \( blockList \) that can fit within the selected space. Then, calculate the score for each block using the formula in \eqref{bscore}. Finally, select the top \( w \) blocks based on their scores, prioritising higher-scoring blocks.

\begin{equation}
	blockScore = \frac{blockArea}{spaceArea} + \frac{b}{avgHigh} \label{bscore}
\end{equation}

In formula \eqref{bscore}, \( blockArea \) denotes the area of the block, \( spaceArea \) represents the area of the selected space, \( avgHigh \) refers to the average length of the remaining boxes after placing the block and \(b\) is a parameter.

\emph{Generate new states proceeds:} Remove the selected space from \( spaceList \), add the newly generated spaces resulting from block placement, and update all spaces that overlap with the block. Remove blocks from \( blockList \) that can no longer be formed with the remaining boxes, and update \( avaibox \) accordingly.

\emph{Greedy approach proceeds:} A space is selected from \( spaceList \), and a block is chosen from \( blockList \). The selected block is placed into the chosen space. This process is repeated until no spaces remain. Once completed, if the filling area of the current configuration surpasses that of \( bestState_a \), then \( bestState_a \) is updated to reflect the new best solution.

\begin{algorithm}[ht]
	\caption{expand}
	\label{expand}
	\KwIn{$state$, $w$}
	\KwOut{$succ$}
	
	$blocks \gets \emptyset$ \;
	
	\While{$blocks = \emptyset$ \textbf{and} $spaceList \neq \emptyset$}{
		$space \gets$ select a space from $spaceList$ \;
		$blocks \gets$ selet $w$ blocks from $blockList$ \;
		
		\If{$blocks$ is empty}{
			Remove $space$ from $spaceList$ \;
		}
		\Else{
			\textbf{break}
		}
	}
	
	$succ \gets \emptyset$ \;
	
	\ForEach{$block \in blocks$}{
		Place the $block$ into $space$ to generate a new state, and then add the new state into $succ$ \;
	}
	
	\Return{$succ$} \;
\end{algorithm}

\subsection{Multiple container filling}	
After completing the minimum container filling process, it is checked whether all boxes have been placed into container \( C_a \). If all boxes are successfully packed, \( bestState_a \) is considered the optimal solution. If not, the remaining boxes are placed into a new container \( C_b \) with width \( W \) and infinite length using the beam search algorithm in a process called quick filling. This beam search is similar to the one in the minimum container filling process, with the key difference being that the search terminates when the search width \( w \) exceeds the number of remaining boxes. Additionally, each iteration selects the longest remaining box and the lowest, leftmost available space.

The best solution from the quick filling process is recorded as \( bestState_b \). After quick filling, the maximum length \( L_{max} \) is computed, and multiple container lengths \( L_i \) are generated, ranging from \( S/W \) to \( L_{max} \) with a step size of 1. For each generated length \( L_i \), a new container \( C_i \) with dimensions \( L_i \times W \) is created. The beam search algorithm then fills each container \( C_i \). Once completed, the configuration with the shortest filling length is selected as the optimal solution.

The beam search algorithm for each \( C_i \) follows the same procedure as in the minimum container filling process. However, since this method does not guarantee that all boxes are placed, adjustments are made to the greedy approach. Specifically, let the state without remaining space be denoted as \( state_c \). Before updating the optimal solution, it is checked whether \( state_c \) satisfies the condition that all boxes have been loaded. If this condition is met, \( state_c \) is compared with the current optimal solution, and the configuration with the shorter filling length is selected as optimal.

If not all boxes have been placed in \( state_c \), a new container with width \( W \) and infinite length is generated, and the quick filling process is applied to pack the remaining boxes. The optimal solution from this process is recorded as \( bestState_f \). Finally, the solutions from \( bestState_f \) and \( state_c \) are combined vertically to form the final solution. This solution is then compared with the current optimal solution, and the configuration with the shortest filling length is chosen as the optimal result. The details of greedy in this step are shown in Algorithm \ref{greedy}.  

\begin{algorithm}[ht]
	\SetAlgoLined
	\KwIn{$state_c$, $bestState$}
	
	\While{$state_c.spaceList$ is not empty}{
		$space \gets$ select a space from $spaceList$\;
		$block \gets$ select a block from $blockList$\;
		
		\eIf{$block$ is not empty}{
			Place $block$ into $space$ and renew $state$\;
		}{
			remove $space$ from $spaceList$ \;
		}
	}
	
	\eIf{all boxes placed in $state_c$}{
		Compare $state_c$ and $bestState$, and store the one with shorter filling length in $bestState$ \;
	}{
		$state_f$ = do quick filling \;
		Combine \( state_c \) and \( state_f \), compare the result with \( bestState \), and update \( bestState \) with the one that has the shorter filling length \;
	}

	\caption{greedy}
	\label{greedy}
\end{algorithm}

\section{Experiments}
This section presents the experimental setup and results used to validate the effectiveness of the proposed algorithm BSPA. First, we describe the computational environment and the standard datasets used to evaluate the algorithm's performance. Next, we explain the procedure for determining the parameter \( b \) in the block scoring formula \eqref{bscore}. We then analyze the algorithm's running time. Finally, we provide a comparative analysis of the performance of BSPA in both the OF and RF scenarios.

\subsection{Experimental environment and datasets}
The algorithm was executed in the following computational environment: AMD EPYC 7H12 2.6 GHz CPU, 120 GB of RAM, the Ubuntu 22.04 operating system, and Java 8 programming language.

The datasets used in this research are as follows:
\begin{itemize}
	\item C: 21 instances proposed by \cite{Hopper2001}.
	\item N: 13 instances proposed by \cite{Burke2004a}.
	\item NT-N, NT-T: 70 instances generated by \cite{hopper2000two}.
	\item KR: 12 instances proposed by \cite{Kroeger1995}.
	\item BWMV: 500 instances proposed by \cite{Berkey1987}, and \cite{Martello1998}.
\end{itemize}

\subsection{Parameter testing}
\begin{figure}[h]
	\centering
	\includegraphics[width=\linewidth]{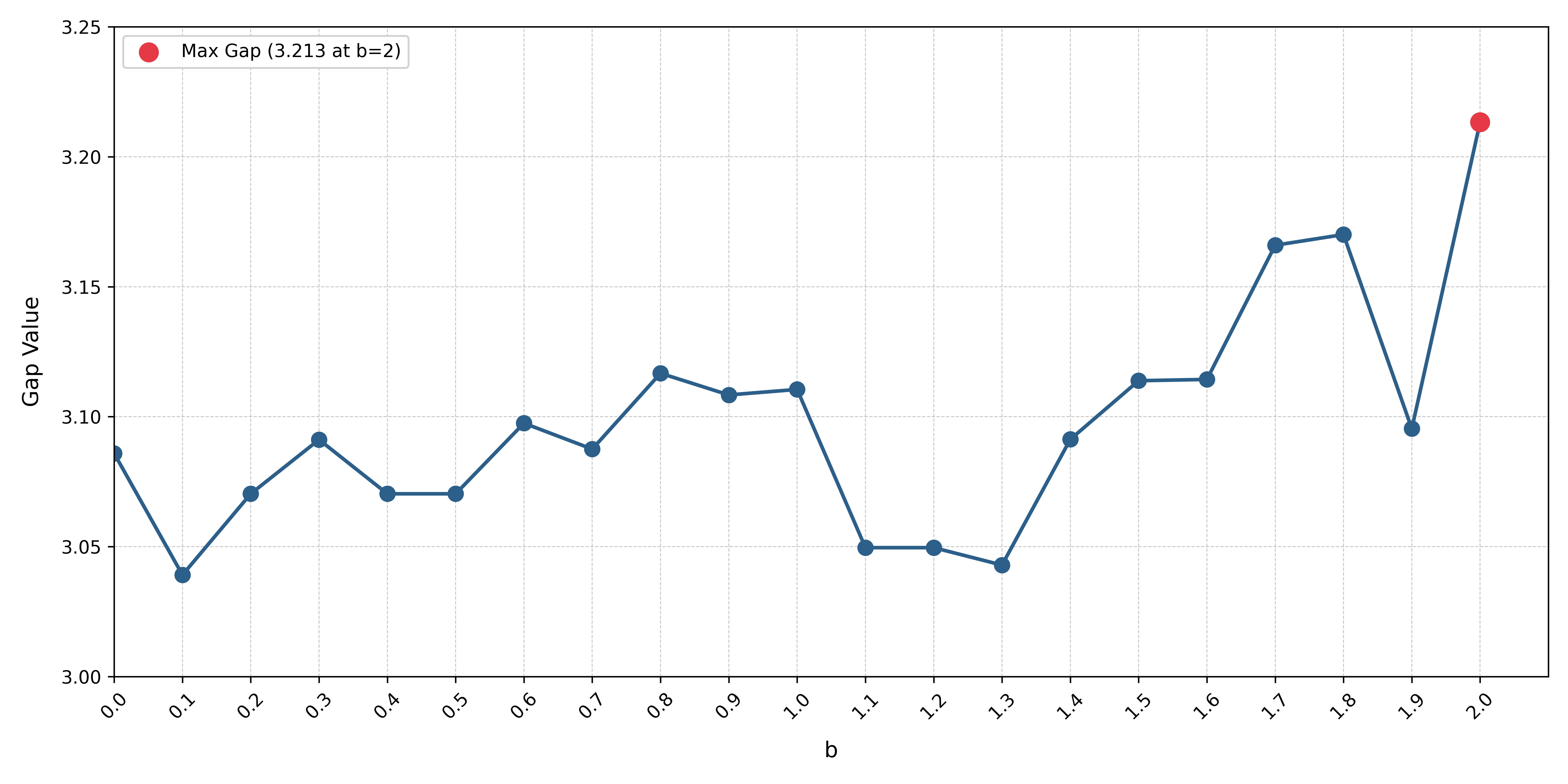}
	\caption{Gap Value Variation with Parameter \(b\).}
	\label{parab}
\end{figure}

To determine the optimal value of the parameter \( b \), 10 problems were randomly selected from each of the aforementioned datasets for testing. Figure \ref{parab} shows the average gap values for different \( b \) values ranging from 0 to 2.0. Notably, the point where \( b = 0 \) indicates that the block's score is solely determined by \( blockArea/spaceArea \). A lower gap value corresponds to a shorter container filling length. The gap is calculated using the formula in \eqref{gap}, where \( L_f \) represents the filling length of the container.

\begin{equation}
	gap = (\frac{L_f}{\frac{S}{W}} -1) \times 100   \label{gap}
\end{equation}

Figure \ref{parab} shows that the gap values are relatively low when \( b \) is set to 0.1, 1.1, 1.2, and 1.3. To further investigate the optimal value of \( b \), we conducted additional tests using these various \( b \) values on the standard dataset. Regarding search time, the beam search algorithm was constrained to a runtime of 30 seconds in both the minimum container filling and multiple container filling processes. The test results are summarized in Table \ref{tab1}.

\begin{table}[htbp]
	\caption{Gap Value Variation with Parameter \(b\).}
	\label{tab1}
	\begin{tabular*}{\tblwidth}{@{}LCCCC@{}}
		\toprule
		\multirow{2}{*}{Datasets} & \multicolumn{4}{c}{$b$} \\
		\cmidrule{2-5}
		& 0.1 & 1.1 & 1.2 & 1.3 \\ 
		\cmidrule{1-5}
		C & \textbf{1.395}  & 1.375  & 1.494  & 1.395   \\ 
		N & \textbf{0.683}  & 0.708  & 0.708  & 0.708   \\ 
		NT-N & \textbf{2.657}  & 2.714  & 2.686  & 2.700   \\ 
		NT-T & \textbf{2.657}  & 2.786  & 2.757  & 2.743   \\ 
		Avg. & \textbf{2.156}  & 2.217  & 2.222  & 2.202   \\ 
		\midrule
		KR & \textbf{2.704}  & 2.704  & 2.619  & 2.556   \\
		BWMV & \textbf{7.623}  & 7.638  & 7.634  & 7.631   \\ 
		Avg. & \textbf{7.508}  & 7.523  & 7.517  & 7.512   \\ 
		\bottomrule
	\end{tabular*}
\end{table}
From Table \ref{tab1}, it can be observed that when \( b = 0.1 \), the algorithm consistently produces a lower gap value across all datasets compared to other values of \( b \). Therefore, \( b = 0.1 \) is chosen as the optimal value.

\subsection{Runtime analysis}
In this section, we analyze the running time of the algorithm presented in this paper. Upon the generation of the block, the algorithm first executes the minimum container filling process. Let the search time for the beam search in this process be denoted as \( T_1 \). If not all boxes are placed into the container, the quick filling process is triggered, with its running time denoted as \( T_2 \). Subsequently, the multiple container filling process is performed. Let the search time for each container in this phase be denoted as \( T_3 \), the number of containers as \( n \), and the number of parallel searches as \( p \). The total running time of the program can be expressed as:

\begin{displaymath}
	T_1 + T_2 + T_3 \times \lceil \frac{n}{p} \rceil
\end{displaymath}

The quick filling process is typically completed very quickly, so its runtime \( T_2 \) can be considered negligible. Therefore, the overall running time of the program is approximately:

\begin{displaymath}
	T_1 + T_3 \times \lceil \frac{n}{p} \rceil
\end{displaymath}

\begin{figure}[htbp]
	\centering
	\includegraphics[width=\linewidth]{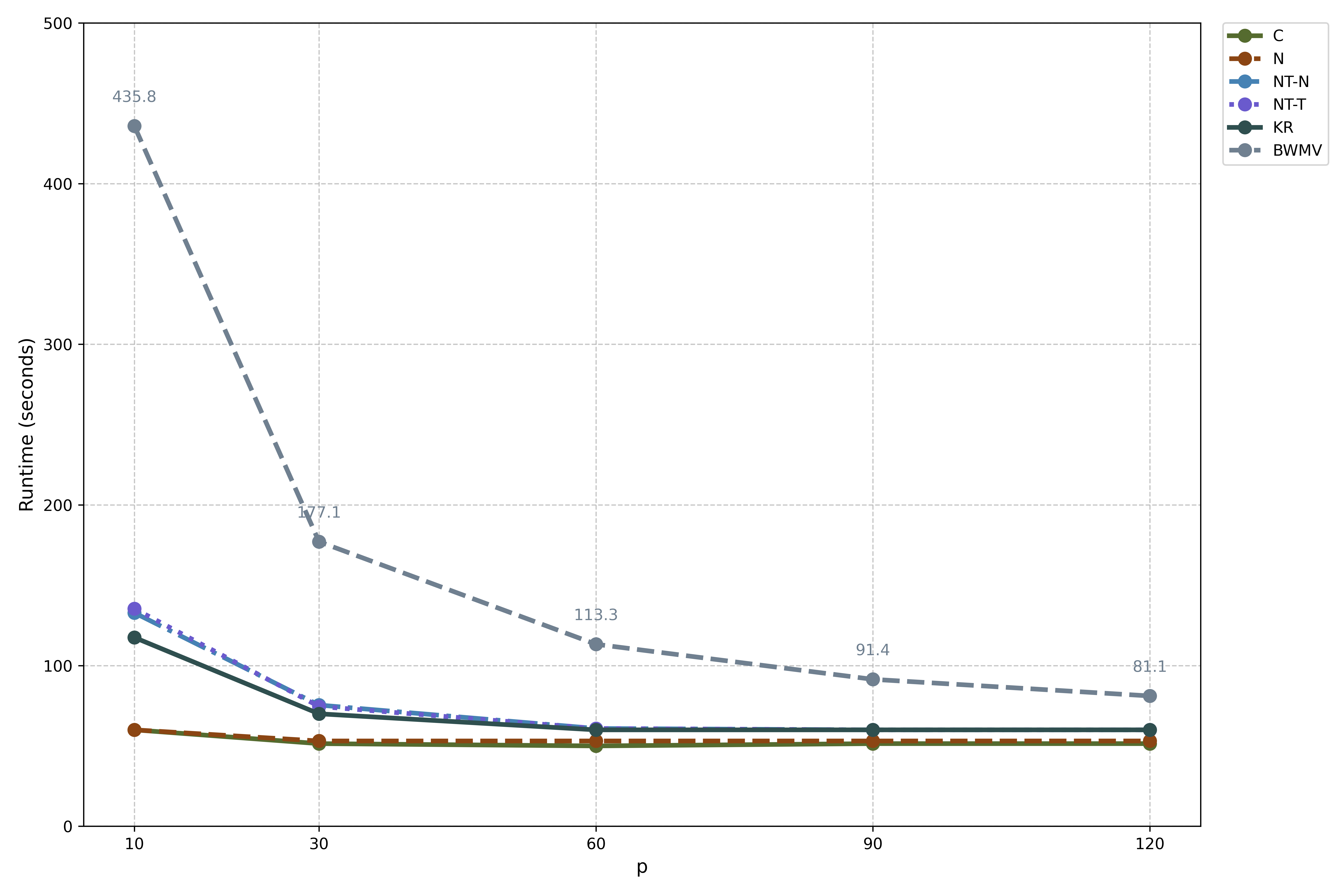}
	\caption{Runtime Variation with Parameter \(p\).}
	\label{pfig}
\end{figure}

It is evident that when \( p \) is sufficiently large, the running time of the entire program approaches \( T_1 + T_3 \). To further investigate the impact of parallel searches on runtime, we set \( T_1 = T_2 = 30 \, \text{s} \) and conducted experiments with various values of \( p \) (specifically 10, 30, 60, 90, and 120) in the OF scenario. The running time for different values of \( p \) is presented in Table \ref{ptable} and Figure \ref{pfig}. As shown in the table, as \( p \) increases, the overall running time of the algorithm gradually approaches \( T_1 + T_3 \) and may even fall below this value. This occurs because, for specific problems, the optimal solution is found during the minimum container filling process, thus preventing the execution of the multiple container filling process and resulting in \( T_3 = 0 \).

However, some counterintuitive results are observed, such as varying gap values for the same \( T_1 \), \( T_3 \), and \( b \) under different values of \( p \). These discrepancies can be attributed to several factors, including the operating system's process scheduling, JVM optimizations, and fluctuations in CPU temperature, all of which can affect computational speed. Therefore, the increase in parallel processes does not have a direct or deterministic relationship with reducing the gap values presented in the table. Instead, the number of parallel processes mainly influences the algorithm's overall running time.
\begin{table*}[htbp]
	\caption{Running Time of BSPA for Different \(p\) Values.}
	\label{ptable}
	\begin{tabular*}{\tblwidth}{@{}LCCCCCCCCCC@{}}
		\toprule
		\multirow{3}{*}{Datasets} & \multicolumn{10}{c}{\(p\)} \\
		\cmidrule{2-11}
		& \multicolumn{2}{c}{10} & \multicolumn{2}{c}{30} & \multicolumn{2}{c}{60} & \multicolumn{2}{c}{90} & \multicolumn{2}{c}{120}\\
		& t(s) & gap & t(s) & gap & t(s) & gap & t(s) & gap & t(s) & gap \\
		\midrule
		C & 60.00  & 1.40  & 51.43  & 1.32  & 50.00  & 1.32  & 51.43  & 1.40  & \textbf{51.43}  & 1.40   \\ 
		N & 60.00  & 0.68  & 53.08  & 0.68  & 53.08  & 0.68  & 53.08  & 0.68  & \textbf{53.08}  & 0.68   \\ 
		NT-N & 132.86  & 2.66  & 75.43  & 2.66  & 60.86  & 2.66  & 60.00  & 2.66  & \textbf{60.00}  & 2.66   \\ 
		NT-T & 135.43  & 2.66  & 74.57  & 2.66  & 60.86  & 2.69  & 60.00  & 2.66  & \textbf{60.00}  & 2.66   \\ 
		Avg & 109.90  & 2.16  & 67.50  & 2.14  & 57.69  & 2.15  & 57.40  & 2.16  & \textbf{57.40}  & 2.16   \\ 
		\midrule
		KR & 117.50  & 2.70  & 70.00  & 2.70  & 60.00  & 2.70  & 60.00  & 2.70  & \textbf{60.00}  & 2.70   \\ 
		BWMV & 435.84  & 7.62  & 177.12  & 7.65  & 113.28  & 7.66  & 91.44  & 7.66  & \textbf{81.12}  & 7.66   \\ 
		Avg & 428.38  & 7.51  & 174.61  & 7.53  & 112.03  & 7.55  & 90.70  & 7.54  & \textbf{80.63}  & 7.54   \\ 			
		\bottomrule
	\end{tabular*}
\end{table*}

\begin{table*}[htbp]
	\caption{BSPA Results Under the OF and RF Scenario.}
	\label{outtable}
	\begin{tabular*}{\tblwidth}{@{}LCCCCCCCC@{}}
		\toprule
		\multirow{3}{*}{Datasets} & \multicolumn{4}{c}{OF} & \multicolumn{4}{c}{RF}\\
		\cmidrule{2-9}
		& \multicolumn{2}{c}{BSPA$_{(0.1,30,30,90)}$} & \multicolumn{2}{c}{BSPA$_{(0.1,60,60,90)}$} & \multicolumn{2}{c}{BSPA$_{(0.1,30,30,90)}$} & \multicolumn{2}{c}{BSPA$_{(0.1,60,60,90)}$} \\
		& gap &  t(s) & gap &  t(s) & gap &  t(s) & gap &  t(s) \\
		\midrule 
		        C & 1.40  & 51.43  & 1.32 & 100.00 & 0.79  & 45.7  & 0.75  & 91.43   \\ 
		N & 0.68  & 53.08  & 0.67 & 106.15 & 0.72  & 53.1  & 0.60  & 101.54   \\ 
		NT-N & 2.66  & 60.00  & 2.59 & 120.00 & 2.14  & 60.0  & 2.04  & 120.00   \\ 
		NT-T & 2.66  & 60.00  & 2.63 & 120.00 & 2.13  & 60.0  & 2.06  & 120.00   \\ 
		Avg. & 2.16  & 57.40  & 2.10 & 114.23 & 1.69  & 56.3  & 1.61  & 111.92   \\ 
		\midrule 
		KR & 2.70  & 60.00  & 2.67 & 120.00 & 1.60  & 60.0  & 1.60  & 120.00   \\ 
		BWMV & 7.62  & 91.44  & 7.56 & 181.68 & 3.03  & 72.7  & 2.79  & 141.12   \\ 
		Avg. & 7.51  & 90.70  & 7.44 & 180.23 & 3.00  & 72.4  & 2.76  & 140.63   \\ 
		\bottomrule
	\end{tabular*}
\end{table*}

\subsection{Experiment and analysis}

Table \ref{outtable} presents the computational results of the algorithm under both OF and RF scenarios. In this table, the notation \(\text{BSPA}_{0.1,30,30,90}\) denotes the execution of the algorithm with parameters \(b = 0.1\), \(T_1 = 30 \, \text{s}\), \(T_3 = 30 \, \text{s}\), and \(p = 90\). This notation is consistently applied across other tables.  

As shown in Table 3, under both OF and RF scenarios, the proposed algorithm demonstrates a consistent reduction in the gap value across all datasets as search time increases. However, despite a twofold increase in search time, the reduction in the gap value remains marginal. Nevertheless, the gap values obtained by our algorithm across all datasets are within an acceptable range. Specifically, for the zero-waste datasets C, N, NT-N, and NT-T, the average fill length exceeds the minimum length by only 2.1\%. For the non-zero-waste dataset KR, this deviation is merely 2.67\%, while for the large-scale dataset BWMV, it remains limited to 7.57\%. These results substantially outperform those of manual filling, indicating that the proposed method can be effectively applied to industrial material cutting to assist enterprises in reducing material waste during the cutting process.

In the previous experiments, we used a platform with a 256-core processor and 120 GB of RAM to run our algorithm. However, this does not imply that the algorithm proposed in this paper is limited to such high-end configurations, as this setup was chosen solely to meet the specific requirements of the experiments. Our algorithm is capable of running on personal computers used in daily life and can even outperform the results obtained in the aforementioned experiments because personal computers usually have higher CPU frequency. Furthermore, if the algorithm is implemented using more efficient programming languages, such as C or C++, the performance of loading results within the same time limit will be further enhanced. Table \ref{persontable} presents the algorithm's performance on a widely used PC with an Intel i7-13700K CPU, 16GB of RAM, and the Windows 11 operating system, using standard datasets.

\section{Conclusions}
This paper presents a beam search-based parallel method for solving the two-dimensional strip packing problem. Experimental results demonstrate that the proposed algorithm can find an efficient placement strategy to minimize the used length of the strip as much as possible. To facilitate future research, the code and datasets are available at \url{https://github.com/Yzhjdj/BSPA}

\begin{table}[htbp]
	\caption{Experimental Results of the Personal Computer.}
	\label{persontable}
		\begin{tabular*}{\tblwidth}{@{}LCCCC@{}}
			\toprule
			\multirow{3}{*}{Dataset} & \multicolumn{4}{c}{BSPA$_{(0.1,30,30,16)}$} \\
			\cmidrule{2-5}
			& \multicolumn{2}{c}{OF} & \multicolumn{2}{c}{RF}  \\
			& gap & t(s) & gap & t(s)\\
			
			\midrule
			C & 1.32  & 51.43  & 0.30  & 40   \\ 
			N & 0.53  & 53.08  & 0.65  & 53.08   \\ 
			NT-N & 2.61  & 96.00  & 2.07  & 79.71   \\ 
			NT-T & 2.64  & 100.29  & 2.06  & 75.43   \\ 
			Avg. & 2.10  & 83.08  & 1.53  & 66.92   \\ 
			\midrule
			KR & 2.67  & 95.00  & 1.60  & 62.50   \\ 
			BWMV & 7.57  & 286.44  & 2.77 & 157.08  \\ 
			Avg. & 7.46  & 281.95  & 2.74 & 154.86  \\ 
			
			\bottomrule
		\end{tabular*}
\end{table}

\printcredits

\bibliographystyle{cas-model2-names}
\bibliography{2dsp}

\end{document}